\documentclass{article}

\usepackage[utf8]{inputenc}
\usepackage[]{graphicx}
\usepackage{hyperref}
\usepackage{amsmath}
\usepackage{centernot}
\usepackage{amssymb}
\usepackage[margin=2.5cm,top=2.5cm,bottom=2.5cm]{geometry}
\usepackage{array}

\newcolumntype{L}[1]{>{\raggedright\let\newline\\\arraybackslash\hspace{0pt}}p{#1}}

\newcommand{\ax}{\alpha^{-1}x}

\title{Rigorous bounds on the Hausdorff dimension of Feigenbaum attractors }
\author{Andrew Burbanks, Andrew Osbaldestin, Judi Thurlby }
\date{February 2021}

\makeatletter
\newsavebox\zzz
\def\mystrut{%
\dimen@\wd\zzz
\divide\dimen@\thr@@
\advance\dimen@-\dp\@arstrutbox
\rule\z@\dimen@}
\def\rotatezzz{%
\rotatebox{90}{\rlap{\kern-\dp\@arstrutbox\usebox\zzz}}}
\makeatother

\begin{document}

\maketitle

\section*{Abstract}
We calculate rigorous bounds on the Hausdorff dimension
of the attractor at the accumulation of the period-doubling
cascade for families of maps with quadratic, cubic, and
quartic critical point.
To do this, we express the attractors as the limit sets
of appropriate Iterated Function Systems constructed
using rigorous bounds on the corresponding renormalisation
fixed point functions.
We use interval arithmetic with rigorous directed rounding modes to show  that the respective dimensions lie in subintervals of the intervals $(0.5370,0.5392)$, $(0.6040,0.6091)$, and $(0.6395,0.6474)$.

\section{The Feigenbaum attractor}

\begin{figure}
\begin{center}
\begin{tabular}{cc}
\includegraphics[width = 7cm]{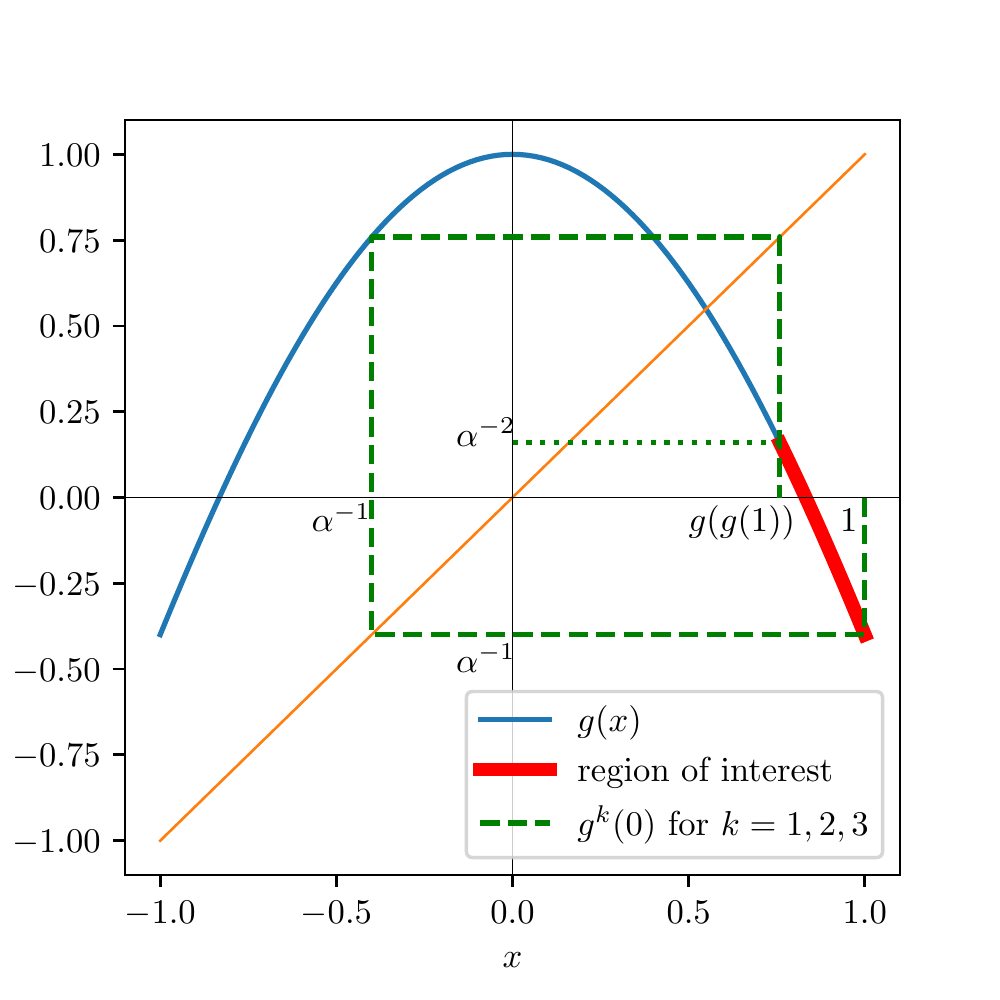} 
     & 
     \includegraphics[width = 7cm]{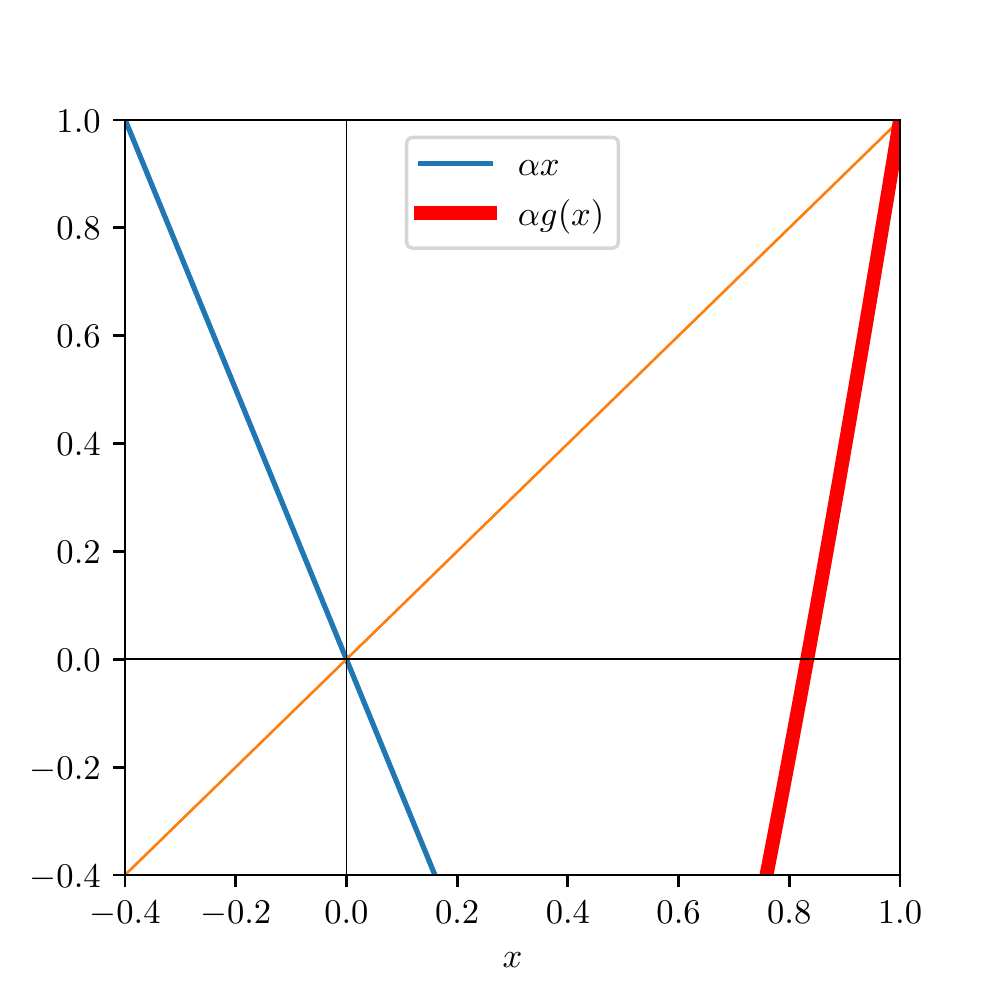}\\ 
    (a) Fixed point function $g(x)$
    & (b) Maps $\alpha x$ and $\alpha g(x)$
 \end{tabular}
\caption{\small (a) Feigenbaum's renormalisation fixed point function (for the case $d=2$) near the origin. The restriction of $g$ to the interval $J=[g(g(1)),1]$ where $g(g(1))\approx 0.76$ is plotted as a thick line (red in colour copy), and
(b) $\alpha x$ (thin line, blue in colour copy) and $\alpha g(x)$ (thick line, red in colour copy) where $\alpha g(x)$ is plotted over $[g(g(1)),1]$. The vertical axis range is $I=[\alpha^{-1},1]$, where $\alpha \approx -2.5$. The inverses of these maps form the IFS. }
\label{FeigenDim}
\end{center}  
\end{figure}

In \cite{Fei88}, Feigenbaum showed that the period-doubling attractor, $A$, can be represented as the limit set of an Iterated Function System (IFS) defined by the two contractive maps: 
\begin{align}
\Psi_0: x &\mapsto \ax,\label{eq:ifs1}\\
\Psi_1: x &\mapsto g^{-1}(\ax),\label{eq:ifs2}
\end{align}
on the interval $I=[\alpha^{-1},1]$,  
where $g$ is the fixed point of Feigenbaum's renormalisation operator \cite{Fei78} 
\begin{equation}
R: g(x)\mapsto \alpha g(g(\alpha^{-1}x)),\quad \text{with } \alpha = g(1)^{-1}.
\label{R}
\end{equation}
Note that here, and in all that follows, $g^{-1}$ denotes the inverse of the  restriction of $g$ to the interval $J=[g(g(1)),1]$, for which $\alpha g(J)=I$ (see figure \ref{FeigenDim}a).

Upper and lower bounds on the Hausdorff dimension of the attractor may be obtained from properties of the constituent maps of the IFS by using  results described in \cite{Fal90} together with rigorous computations of bounds on the renormalisation fixed-point, $g$.
The operator $R$ is defined on suitable sets of functions having a critical point of (even) integer degree $d$ at the origin. Note that the choice of $\alpha$  preserves the normalisation $g(0)=1$.  We note that both $g$ and $\alpha$ depend on $d$; in what follows we suppress this dependence in the notation, except where needed, for clarity.

The initial interval of the IFS is determined by the points along the forward orbit,
\begin{equation}
 \left(g^{k}(0)\right)_{k\ge 0}
= \left(0,1,\alpha^{-1},g(\alpha^{-1}),\alpha^{-2},...\right), 
\end{equation}
of the critical point.
In fact, the endpoints of the subintervals in successive generations of the IFS comprise this forward orbit. The functional form (equations \ref{eq:ifs1},\ref{eq:ifs2}) of the IFS is common across families of maps for any integer degree $d\ge 2$ critical point. 

Estimates for the Hausdorff dimension of the  attractor in the case of maps with  critical point of degree $2$ have been provided by Grassberger \cite{Gra81} \cite{Gra85}, Bensimon, Jensen and Kadanoff \cite{Ben86}, Kovacs \cite{Kov89} (by expressing the dimension as an eigenvalue of a Perron-Frobenius operator), amongst others.
In particular, the cycle expansion technique of Cvitanovic and coworkers \cite{Cvi93} has been extremely profitable culminating in the best numerical result to date by Christiansen et al \cite{Chr90} yielding $27$ digits:
\begin{equation}
    \mathrm{dim_H}(A_2)\approx  0.538\,045\,143\,580\,549\,911\,671\,415\,567,
\end{equation}
where $A_d$ denotes the attractor corresponding to the universality class of maps with degree $d$ critical point.
The approach that we take, making use of the IFS along with rigorous bounds on the relevant renormalisation fixed points, delivers rigorous bounds on the dimension but, due to computational cost (and inevitable looseness of bounds), is unable to match the precision of those numerical techniques.

\begin{figure}
\begin{center}
    \includegraphics[width = 12cm]{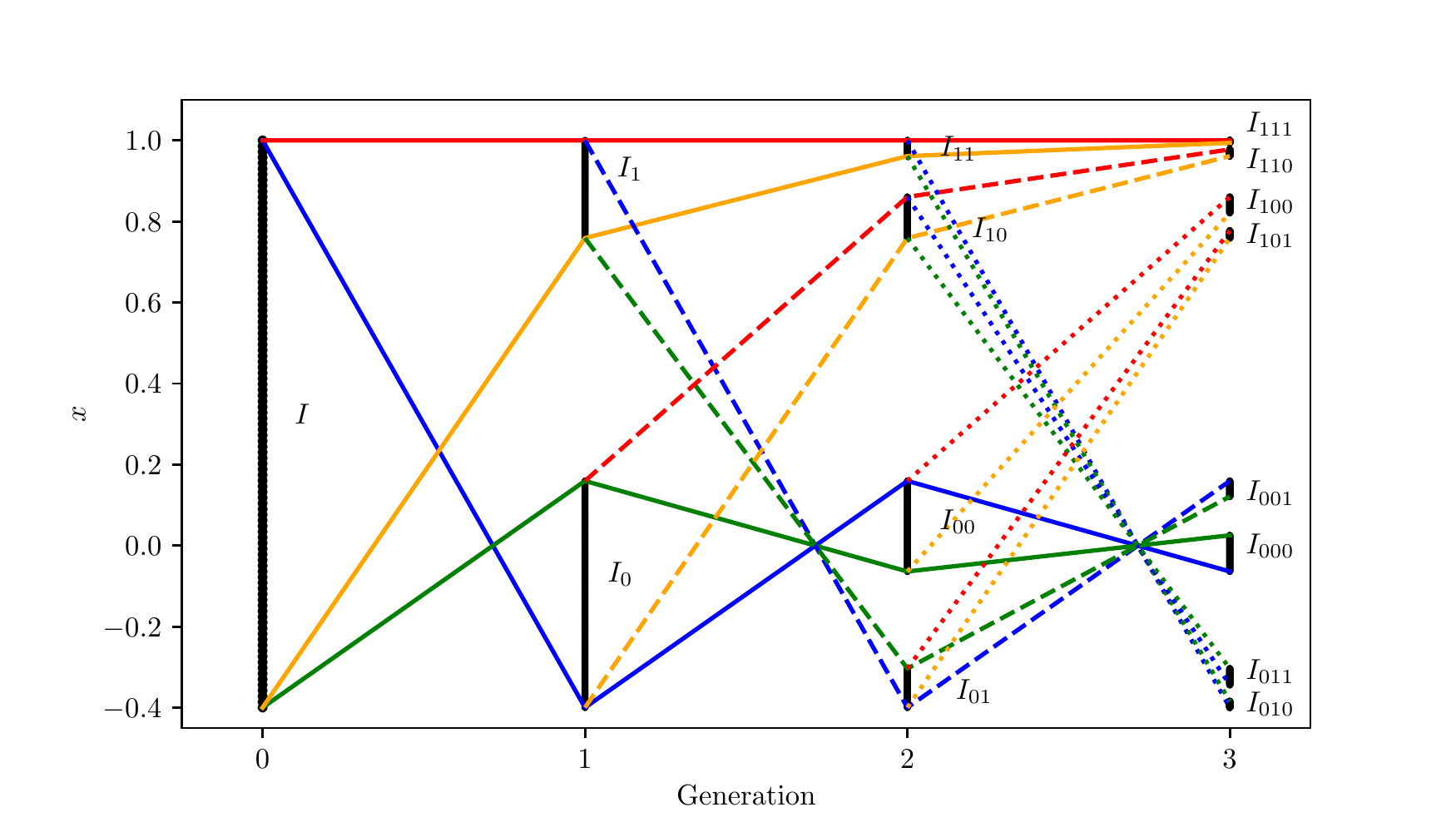} 
    \caption{The first three generations of the IFS in the case  $d=2$, illustrating the iterative construction of the  attractor. 
    To aid understanding we emphasise the parts of the construction corresponding to the symbol sequence $\{000\}$ (blue/green in colour version) and sequence $\{111\}$ (red/orange) with solid lines (see section \ref{rigFeigDim}). The dashed lines indicate a change of symbol after the first generation and the dotted lines indicate a change after the second generation. The map $\Psi_0$ (blue/green), is order-reversing, and $\Psi_1$ (red/orange) is order-preserving.} 
    \label{fig:QuadraticIFS}
    \end{center}
\end{figure}

Figure \ref{fig:QuadraticIFS} shows the subintervals corresponding to the first three generations of the IFS for the Feigenbaum attractor in the case $d=2$.

\section{Hausdorff dimension}\label{rigFeigDim}

Following the method described in \cite{Fal90} we use bounds on the derivatives of the IFS functions to bound the Hausdorff dimension of the attractor.
At generation $n$ of the IFS, we have $2^n$ maps, $\Psi_\sigma$ indexed by symbol sequences $\sigma=a_{n-1}\ldots a_1a_0\in\{0,1\}^n$, as shown in figure \ref{fig:QuadraticIFS} (for $n\le 3$), defined by
\begin{equation}
\Psi_{\sigma}=\Psi_{a_{n-1}\ldots a_1a_0}=
\Psi_{a_{n-1}}
\circ \cdots
\circ \Psi_{a_1}
\circ \Psi_{a_0}.
\end{equation}
The principal ingredient of the IFS approach to calculating dimensions is knowledge of the contractivity and `coercivity' factors of each map.
For each $\Psi_{\sigma}$, uniform contractivity and coercivity constants, $c_{\sigma}$ and $d_{\sigma}$ respectively, satisfy
\begin{equation}
    0<d_\sigma\le\frac{|\Psi_\sigma(x)-\Psi_\sigma(y)|}{|x-y|}\le c_\sigma<1 \qquad \forall x,y\in I
    \ \mbox{with $x\neq y$.}
\end{equation}
We calculate bounds on the derivatives of each map. Smoothness of the $\Psi_\sigma$ gives
\begin{equation}
    \inf_{z\in I}|\Psi_\sigma'(z)|\le
    \frac{|\Psi(x)-\Psi(y)|}{|x-y|}\le
    \sup_{z\in I}|\Psi_\sigma'(z)| \qquad \forall x,y\in I\ \mbox{with $x\neq y$,}
\end{equation}
and thus suitable constants, $c_\sigma$ and $d_\sigma$, may be found by bounding the suprema and infima of the absolute values of the derivatives of the maps.
Solving the corresponding partition function equations for $s_n$ and $r_n$,
\begin{equation}
    \sum_{\sigma}c_\sigma^{s_n} = 1,\quad
    \sum_{\sigma}d_\sigma^{r_n} = 1,
\end{equation}
at each generation $n$ then gives bounds
\begin{equation}
    r_n\le\mathrm{dim_H}(A)\le s_n,
\end{equation}
on the Hausdorff dimension of the attractor, $A$, of the IFS.

The derivatives of the IFS functions are given by:
\begin{align}
\Psi_0': x &\mapsto \alpha^{-1},\\
\Psi_1': x &\mapsto \alpha^{-1}(g^{-1})'(\alpha^{-1}x),
\end{align}
defined on the interval $I=[\alpha^{-1},1]$, with $\alpha = g(1)^{-1}$.

\begin{figure}
\begin{center}
\begin{tabular}{cc}
\includegraphics[width = 5.75cm]{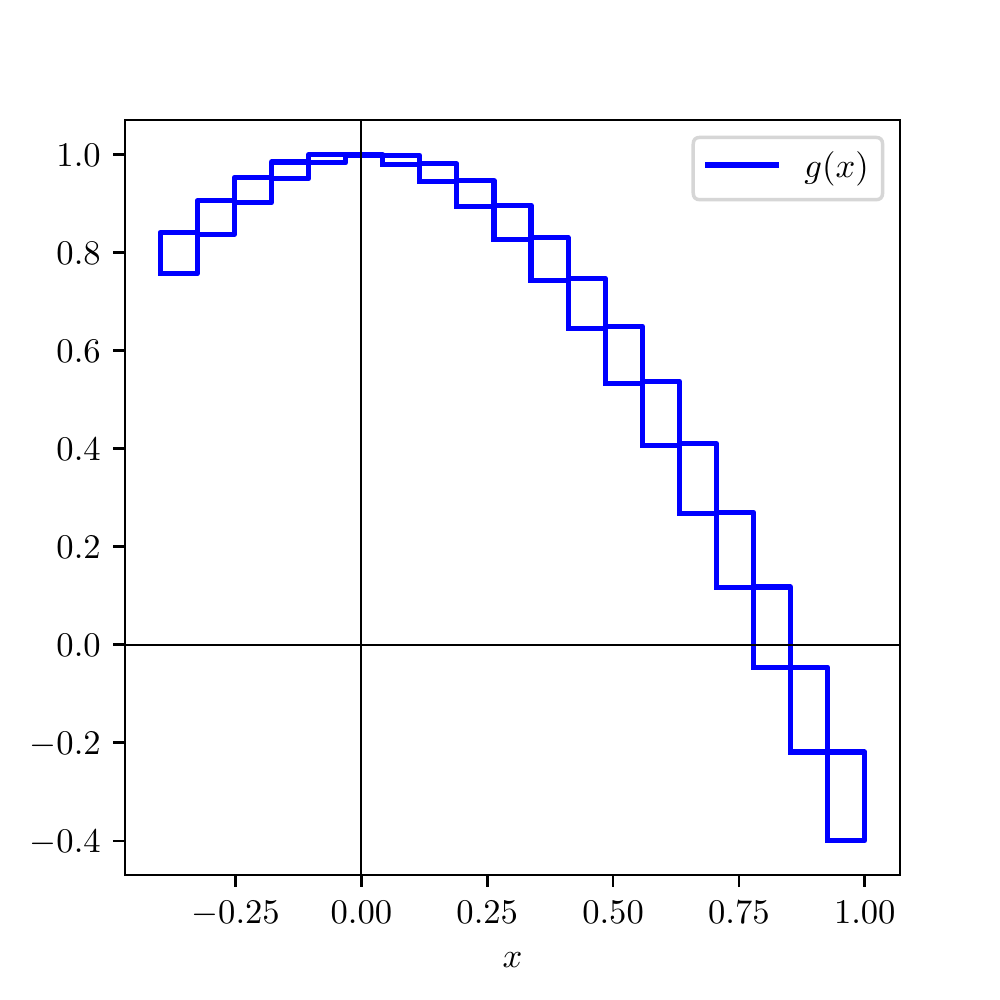} 
& \includegraphics[width = 5.75cm]{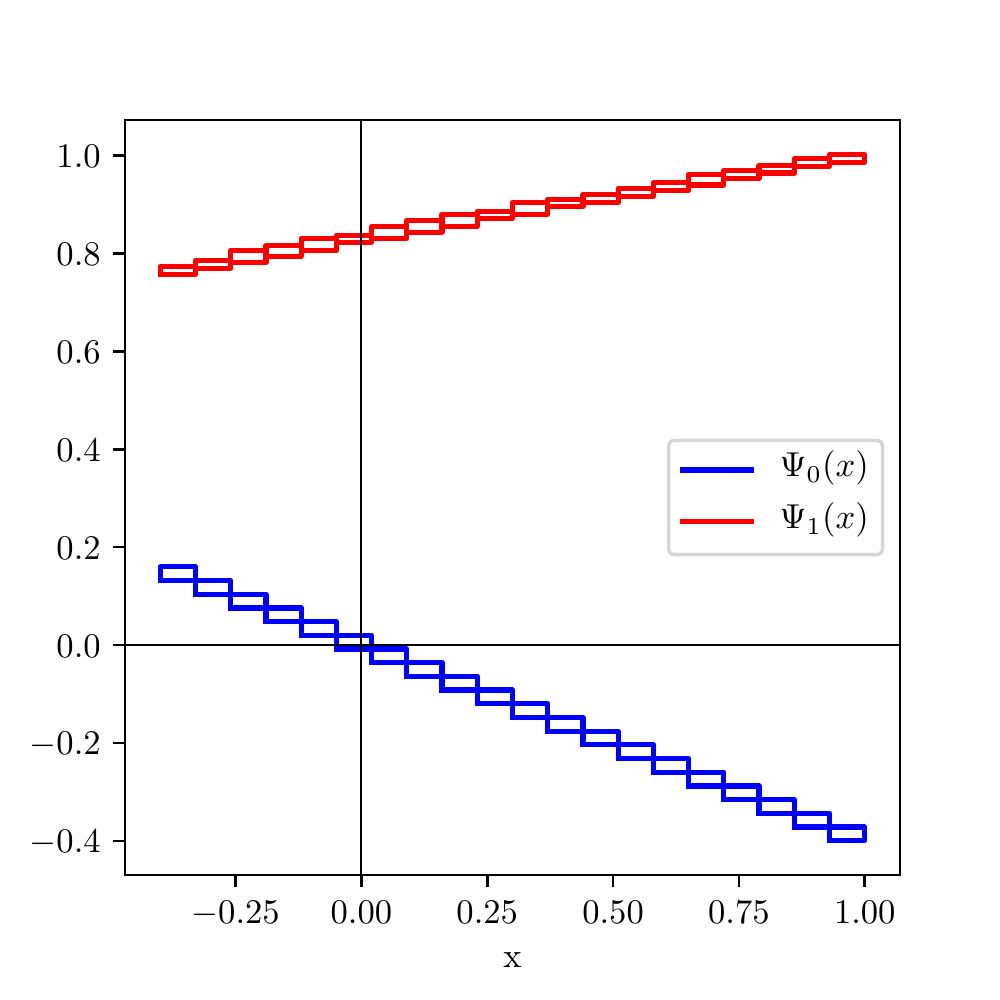}\\ 
(a) Rigorous bounds on $g$
& (d) Rigorous bounds on $\Psi_0$ and $\Psi_1$\\
    
     \includegraphics[width = 5.75cm]{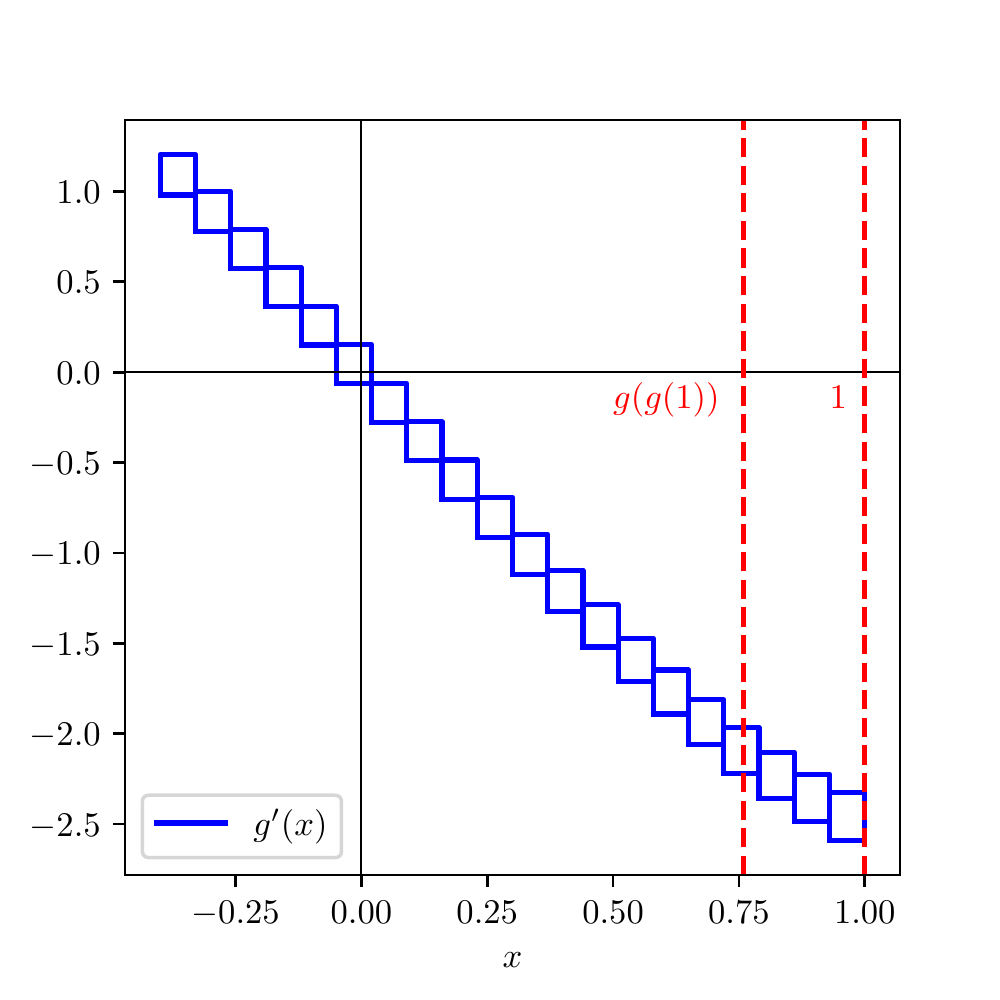}
     & \includegraphics[width = 5.75cm]{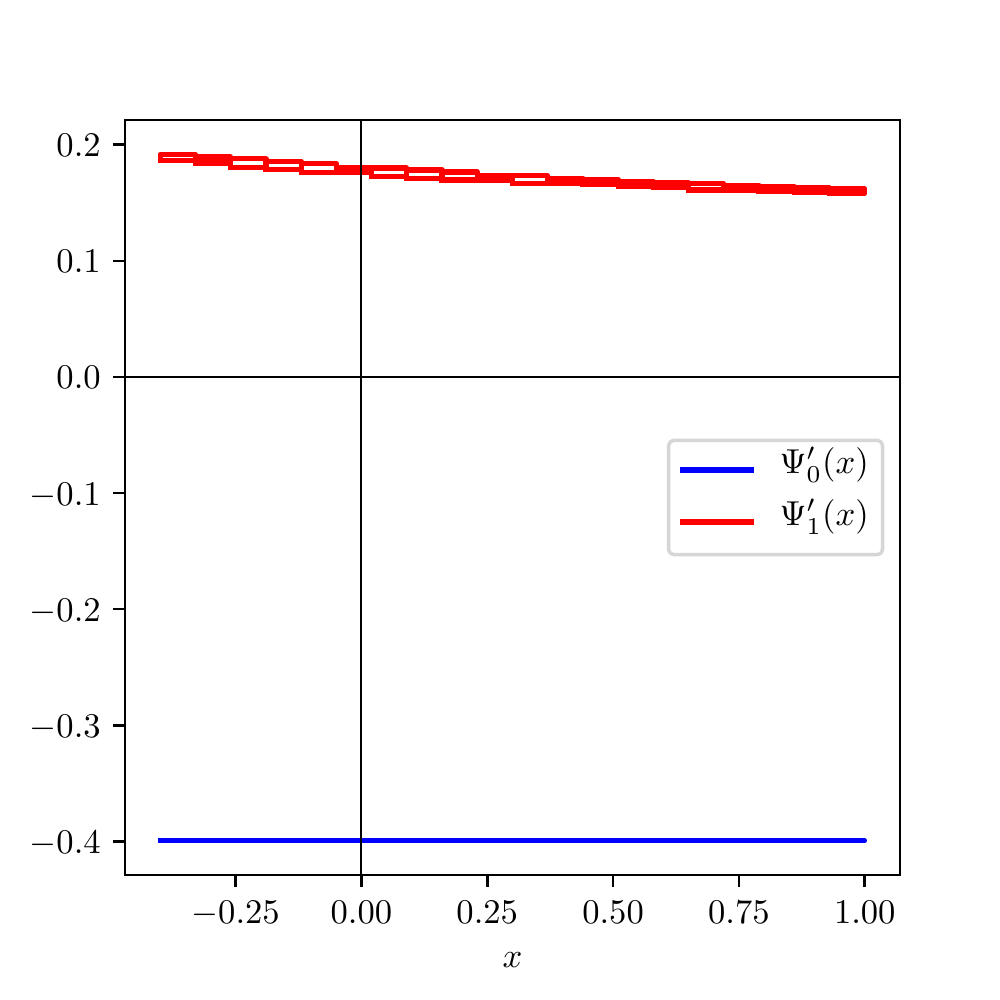}\\ 
    (b) Rigorous bounds on $g'$
    & (e) Rigorous bounds on $\Psi_0'$ and $\Psi_1'$\\
    
    \includegraphics[width = 5.75cm]{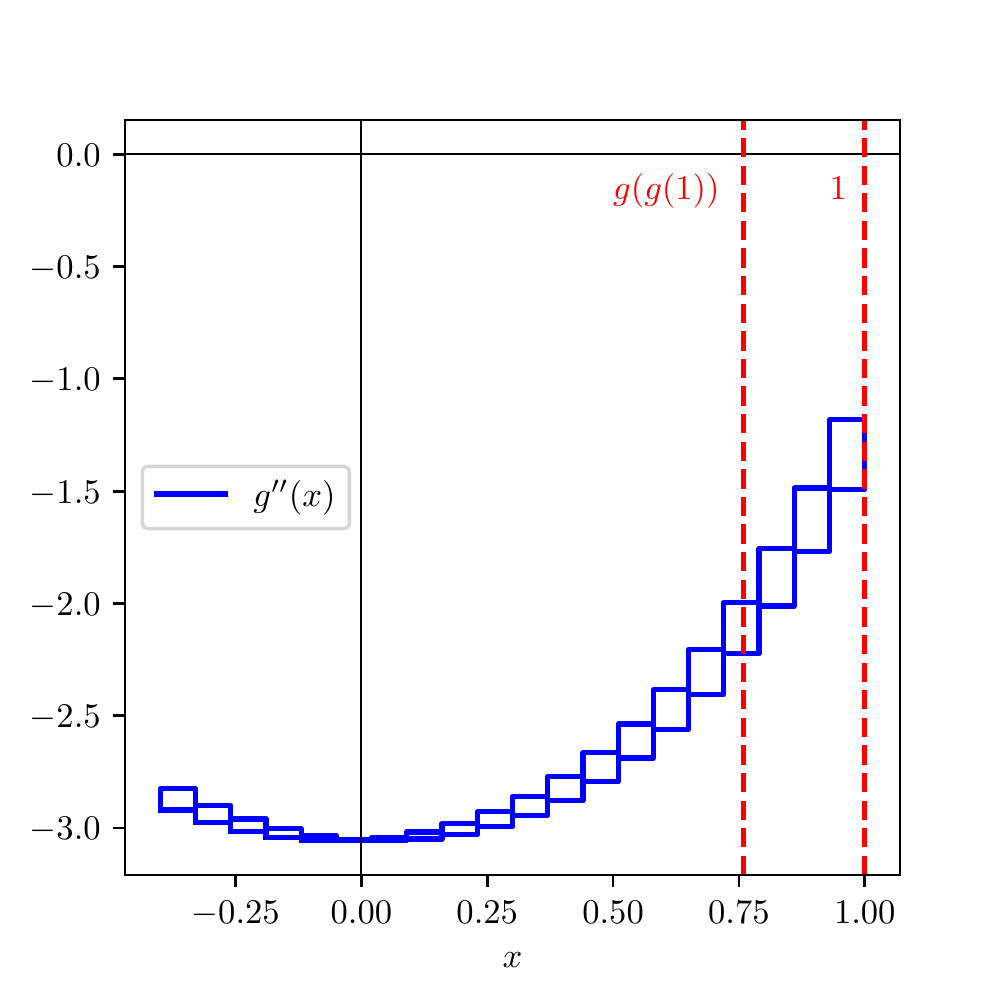} 
    & \includegraphics[width = 5.75cm]{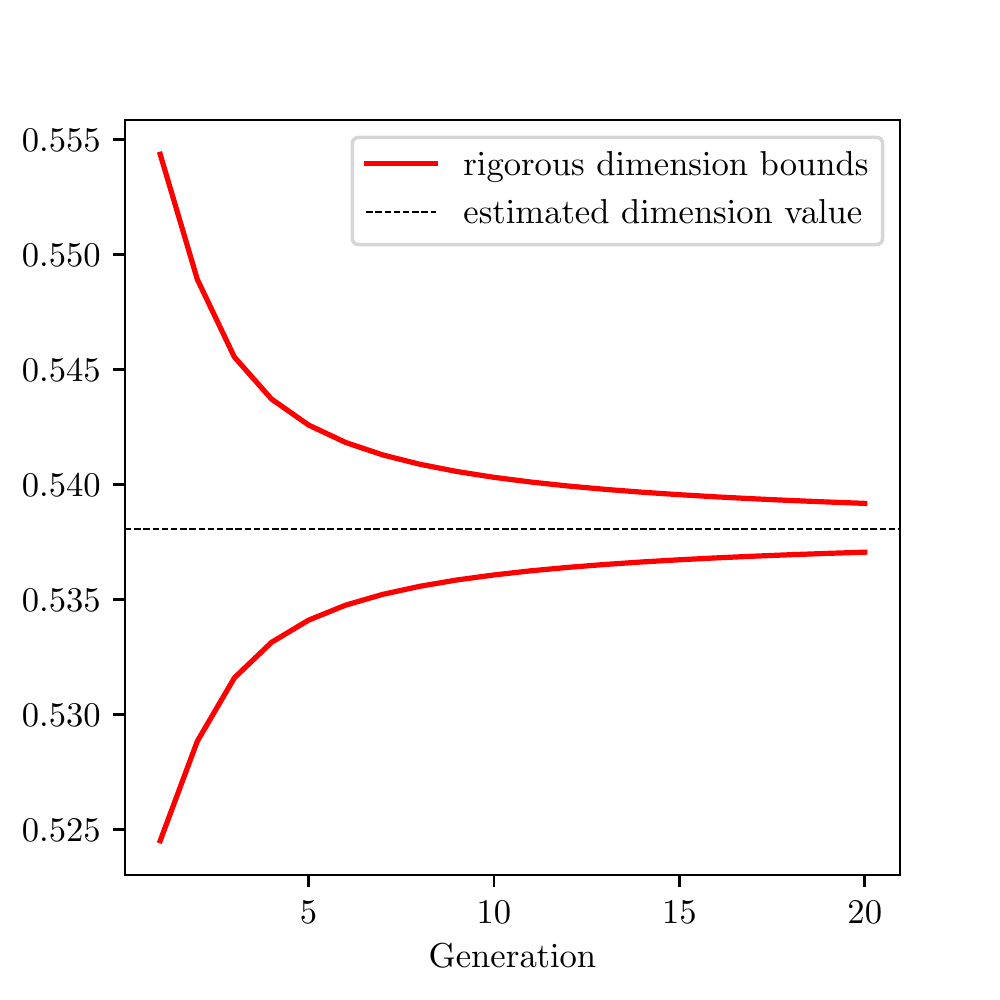}\\ 
    (c) Rigorous bounds on $g''$
    & (f) Rigorous bounds on dimension
   
\end{tabular}
\caption{Construction of the IFS maps and resulting bounds on dimension for the case $d=2$. The rectangles provide a coarse rigorous covering of the graphs of the corresponding functions: (a) the fixed point function, $g$; (b) the first derivative of $g$, with the vertical dashed lines indicating the interval $[g(g(1)),1]$; (c) the second derivative of $g$, demonstrating that $g''(x)<0$ on $[g(g(1)),1]$; (d) rigorous bounds on the functions $\Psi_0$ and $\Psi_1$; and (e) rigorous bounds on the derivatives of $\Psi_0$ and $\Psi_1$. (f) shows convergence of rigorous upper and lower bounds on dimension with increasing IFS generation.}
\label{Q2Dim}
\end{center}  
\end{figure}

In order to form the partition function equations for the IFS we require bounds on  contractivity over the fundamental interval of the IFS.
Clearly $\Psi_0$ is monotonically decreasing, with  $\Psi_0'$  constant. However, $\Psi_1'$ presents some challenges.
We need to bound  $(g^{-1})'$ on the interval $\alpha^{-1}I$. It might at first seem that a  covering by sub-intervals is required, however we are able to establish rigorously for $d = 2,3,4$ that $\Psi_1'$ is monotonic on $I$ (illustrated in figure \ref{Q2Dim} for the case $d=2$) and therefore we need only evaluate it at the endpoints of the IFS interval.
We can therefore bound the contractivities and coercivities by evaluating derivatives only at the endpoints of IFS subintervals.

We prove the monotonicity of $\Psi_1'$ by reference to $g''$ (illustrated in figure \ref{Q2Dim}c for the case $d=2$).

Since the second derivative of $g$ is proven to be strictly negative, $g'$ must be strictly monotonically decreasing on the relevant interval $[g(g(1)), 1]=\Psi_1(I)$ (note also that $g'$ is bounded away from $0$ on this interval, establishing monotonicity of $g$ itself). 
To calculate the derivative of the inverse map we use the standard result that $(g^{-1})'(y) = g'(x)^{-1}$ with $y=g(x)$.  Thus we may calculate $g'(x)^{-1}$ instead of $(g^{-1})'(y)$, provided that we have a means to solve $y=g(x)$ for $x$ in a rigorous manner.
We thus prove that $(g^{-1})'$  is monotonically increasing on the interval $\alpha^{-1}I$ and thus $\Psi_1'$ is monotonically decreasing on $I$. 

We use a ball $B$, in a suitable space of analytic functions, previously proven to contain the renormalisation fixed point, $g$, as described in detail in \cite{Bur20deg2}. The ball $B$ is centred on a high-degree polynomial approximation to $g$.
For the implementation of this method we need to address a technical issue when using bounds on $g$ (in the space of functions) to compute bounds on its inverse $g^{-1}$. For example, given a ball of functions $B$ containing an invertible function $f$, it is not sufficient, in order to bound $f^{-1}$, to find another ball $C$ such that computing bounds on $B\circ C:=\{f\circ h:f\in B,h\in C\}$ results in a ball $D\supseteq B\circ C$ containing the identity function. Looseness of the bounds computed for the composition means that this does not imply that $f^{-1}\in C$ in general. To circumvent this, we instead make use of the function ball $B\ni g$ directly, together with a rigorous root-finding method to bound an interval $X=[a,b]\ni x=g^{-1}(y)$ such that $g(X)\ni y$ for given $y$.

As described earlier we need rigorous bounds on $\Psi_1'$ at the endpoints of IFS intervals, which raises a further obstacle: given a ball $B=B(f,r)$ it is not possible, in general, to find a ball $B'$ such that $B'\supseteq\{h':h\in B\}$, as the corresponding derivative operator is unbounded.
We note, however, that we require bounds valid only for the renormalisation fixed point $g$ itself. Thus, to overcome this problem, we may make use of the fixed point equation,
\begin{align}
    g(x) &=\alpha g(g(\alpha^{-1}x)),
\end{align}
differentiating both sides to give
\begin{align}
    g'(x) &=g'(g(\alpha^{-1}x))g'(\alpha^{-1}x),
\end{align}
and note that the composition of a derivative with a function, $f'\circ g$, is indeed bounded within a ball of functions for suitably-bounded functions $f,g$, as shown in \cite{Eck84}.

When bounding the contractivities and coercivities of the IFS maps in successive generations, we use the chain rule: at generation $n+1$ we have
\begin{align}
    I_{\sigma_{n+1}\sigma}=&(\Psi_{\sigma_{n+1}}\circ \Psi_{\sigma})(I) = \Psi_{\sigma_{n+1}}(I_{\sigma}),\\
    \Psi'_{\sigma_{n+1}\sigma}=&(\Psi_{\sigma_{n+1}}\circ \Psi_{\sigma})' = (\Psi'_{\sigma_{n+1}}\circ \Psi_{\sigma})\cdot \Psi'_{\sigma}.
\end{align}
For instance, at generation $2$ we require the derivatives $\Psi'_{ij}$ evaluated at the endpoints of $I$ for $i,j \in \{0,1\}$. We have 
\begin{equation}
    \Psi'_{ij}=(\Psi_i\circ \Psi_j)'=(\Psi_i'\circ \Psi_j)\cdot \Psi_j',
\end{equation}
where $\Psi_j'$ has already been computed at the endpoints of $I$ at generation $n=1$ and $\Psi_i'$ must then be computed at the endpoints of the subinterval $I_j=\Psi_j(I)$.

The partition function equations are solved using rigorous interval arithmetic to give new upper and lower bounds on the Hausdorff dimension of the attractor at each generation.

Using a ball of functions, $B$, with $\ell_1$ radius $10^{-9}$ centred on a polynomial  of degree 40   that we previously proved contains the fixed point $g$ \cite{Bur20deg2}, using the techniques of \cite{Eck84}, after 20 generations of the IFS we obtain rigorous bounds on the Hausdorff dimension of the Feigenbaum attractor $A_2$: 
\begin{equation}
0.5370523555103920147606639205 < \mathrm{dim_H}(A_2) < 0.5391744736510156113653618025.
\end{equation}
Figure \ref{Q2Dim}f shows the convergence of the upper and lower bounds of the dimension.

\section{Hausdorff dimension - degree 3 critical point}
\label{chD3Dim}

\begin{figure}
\begin{center}
\begin{tabular}{cc}
\includegraphics[width = 5.75cm]{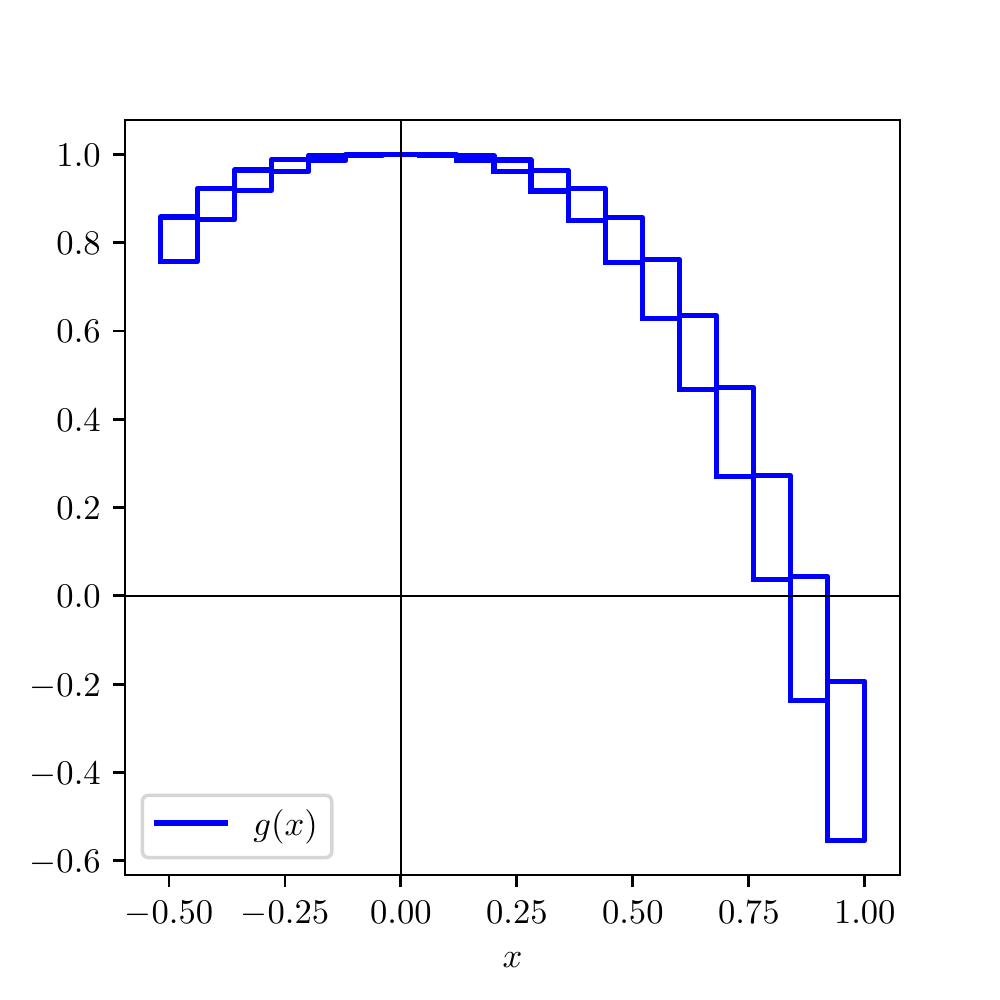}
& \includegraphics[width = 5.75cm]{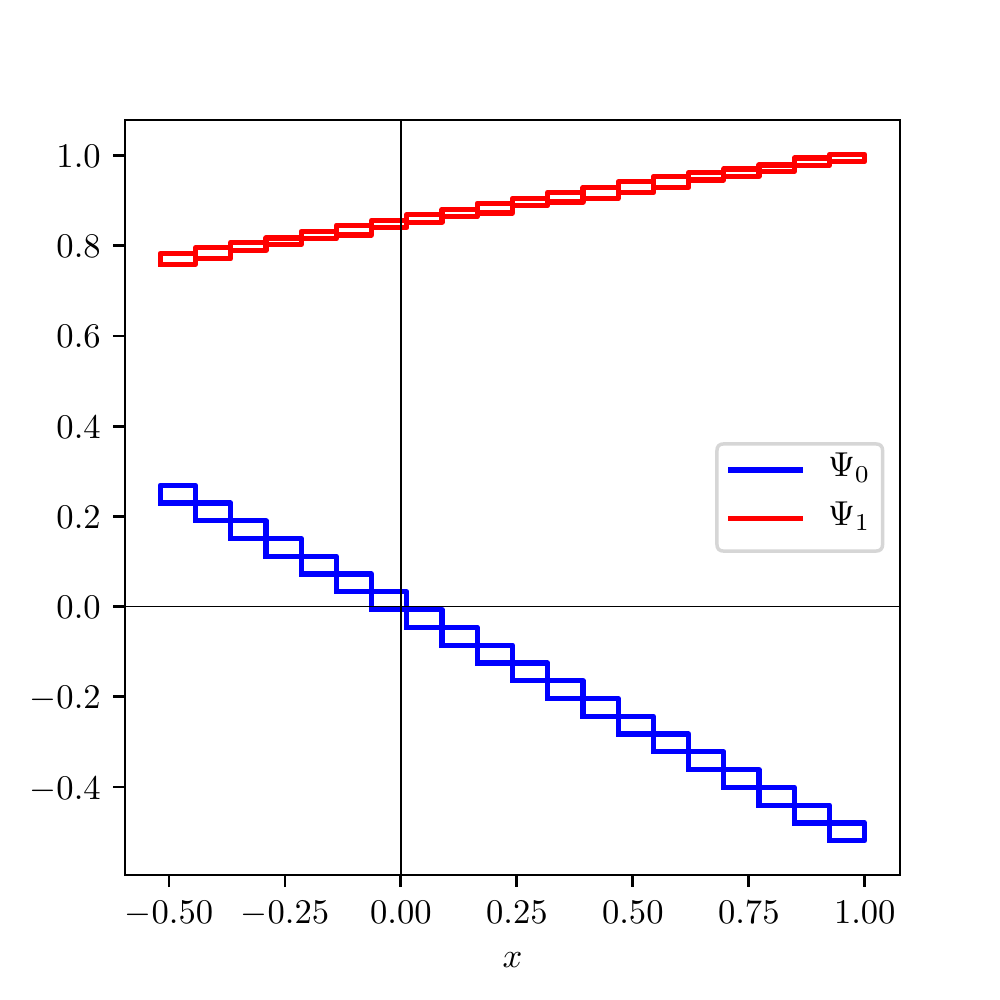}\\
(a) Rigorous bounds on $g$
& (d) Rigorous bounds on $\Psi_0$ and $\Psi_1$\\
    
     \includegraphics[width = 5.75cm]{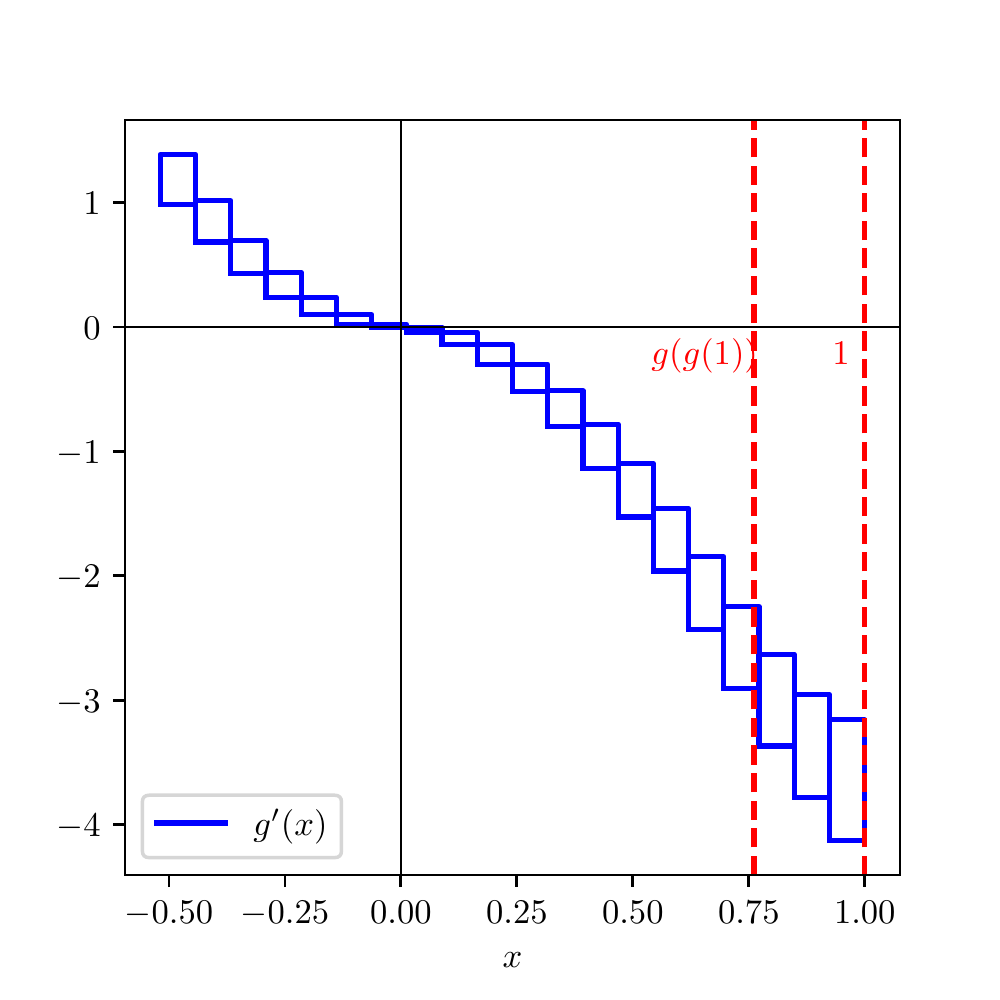}
     & \includegraphics[width = 5.75cm]{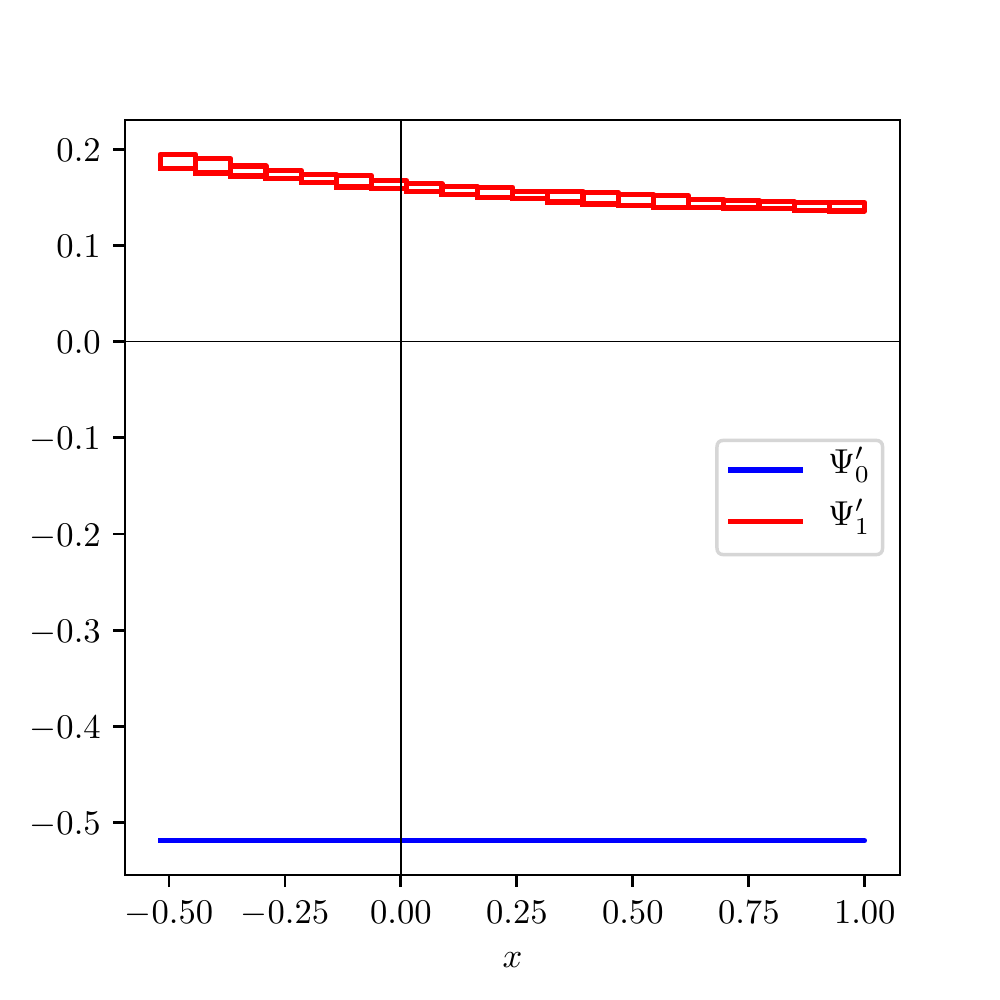}\\
    (b) Rigorous bounds on $g'$
    & (e) Rigorous bounds on $\Psi_0'$ and $\Psi_1'$\\
    
    \includegraphics[width = 5.75cm]{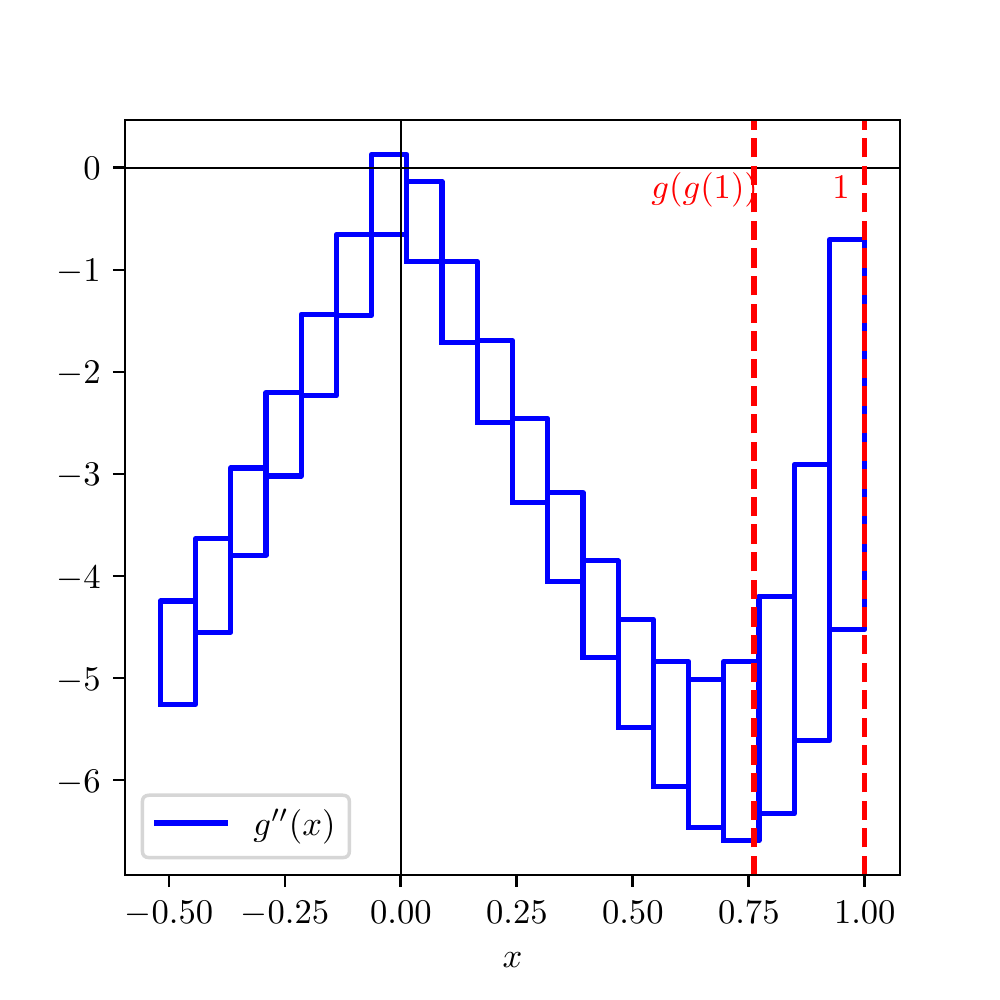}
    & \includegraphics[width = 5.75cm]{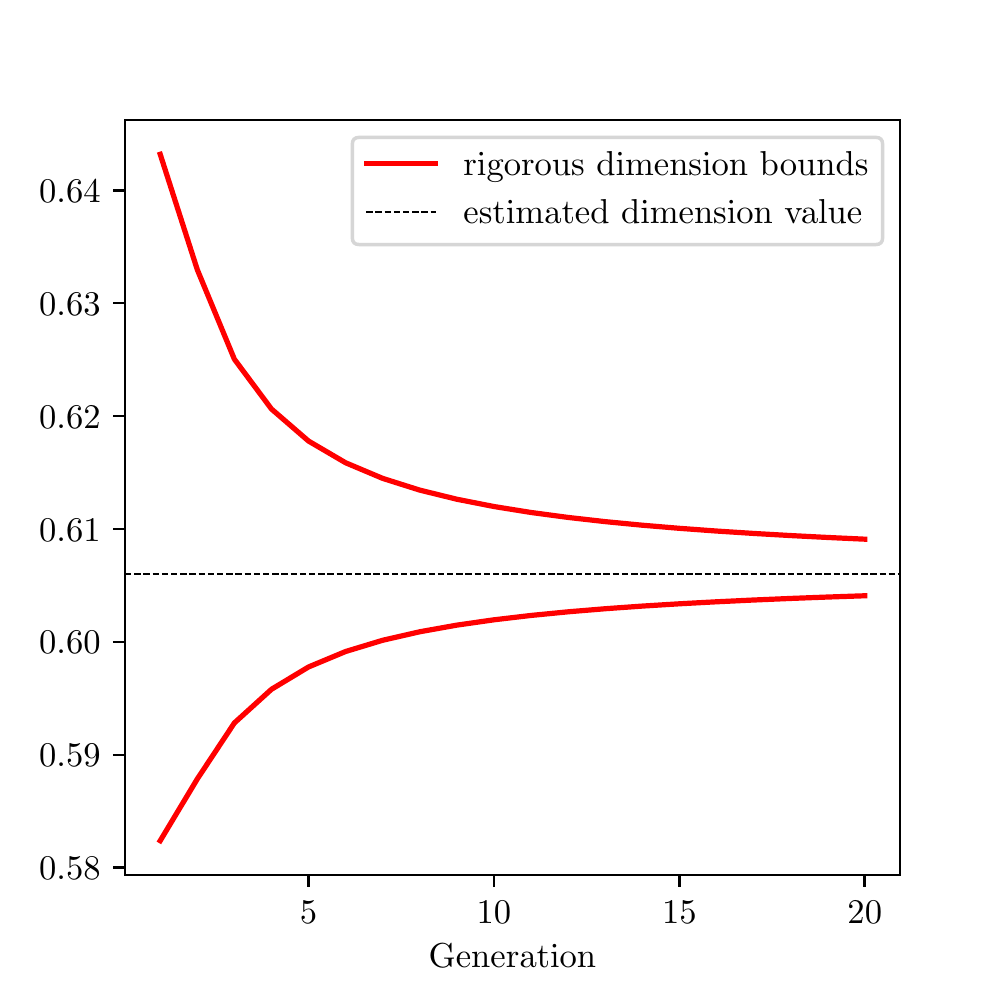}\\
    (c) Rigorous bounds on $g''$
    & (f) Rigorous bounds on dimension
   
\end{tabular}
\caption{Construction of the IFS maps and resulting bounds on dimension for the case $d=3$. The rectangles provide a coarse rigorous covering of the graphs of the corresponding functions: (a) the fixed point function, $g$, for which $g(x)=g_{+}(|x|)$; (b) the first derivative of $g$, with the vertical dashed lines indicating the interval $[g(g(1)),1]$; (c) the second derivative of $g$, demonstrating that $g''(x)<0$ on $[g(g(1)),1]$; (d) rigorous bounds on the functions $\Psi_0$ and $\Psi_1$; and (e) rigorous bounds on the derivatives of $\Psi_0$ and $\Psi_1$. (f) shows convergence of rigorous upper and lower bounds on dimension with increasing IFS generation.}
\label{Q3Dim}
\end{center}  
\end{figure}

We adapt the procedure above to the attractors corresponding to the universality class of maps with odd degree $d$ critical point exemplified by the prototypical family
\begin{align}
x_{n+1}=1-\mu|x_n|^d.
\end{align}
In order to encode the absolute value, we modify the renormalisation operator, equation (\ref{R}), defining separate analytic functions for use with positive and negative operands:
\begin{align}
g(x) =
    \begin{cases}
         g_+(x), & x\geq 0\\
         g_-(x), & x<0,
    \end{cases}
\end{align}
where $g_-(x) = g_+(-x)$ for $x<0$.
The renormalisation operator $R$ becomes:
\begin{align}
R:
\begin{cases}
    g_+(x)\mapsto\alpha g_+(g_-(\alpha^{-1} x)),\\ 
    g_-(x)\mapsto\alpha g_+(g_+(\alpha^{-1} x)),
\end{cases}
\end{align}
where $\alpha=g(1)^{-1}=g_+(1)^{-1}<0.$ 
It suffices to work with just $g_+$ on a carefully chosen domain using a modified operator which with a slight abuse of notation we also refer to as $R$:
\begin{equation}
    R(g_+)(x): g_+(x) \mapsto \alpha g_+(g_+(-\alpha^{-1} x))
    \label{R3}
\end{equation}
where $\alpha=g_+(1)^{-1}$ (note the presence of $-\alpha^{-1}$, rather than $\alpha^{-1}$ in the inner bracket).
The techniques used for $d=2$ above may now be used for maps with odd degree critical point, for which the corresponding map $g(x)=g_+(|x|)$, is not itself analytic at the origin.

The maps defining the corresponding IFS have the same functional form as for $d=2$, equations (\ref{eq:ifs1}) and (\ref{eq:ifs2}), with $\alpha=\alpha_d$ being the relevant universal constant.
We work with a ball of functions proven to contain the corresponding fixed point, $g=g_d$ (in fact, $g_+$), of the renormalisation operator using the method described in detail in \cite{Bur20deg3}.

Figure \ref{Q3Dim} illustrates the results in the case $d=3$. 
We establish rigorously that $g''$ is strictly negative on the domain $[g(g(1)), 1]$ (illustrated in figure \ref{Q3Dim}c), thus proving that $g'$ is strictly monotonically decreasing (note that $g'$ is itself bounded away from zero, establishing monotonicity of $g$).  Thus, $\Psi_1'$ is monotonically decreasing on the domain $[\alpha^{-1},1]$ as required. Again, the significance of this is that we need only bound the derivatives at the endpoints of the intervals in order to bound contractivity of the IFS maps. 

Figure \ref{Q3Dim}f shows the convergence of the rigorous upper and lower bounds of the Hausdorff dimension with each step. The bounds produced verify the previously published numerical estimate for the dimension of the attractor for maps with degree 3 critical point \cite{BHu}:
\begin{equation*}
    \mathrm{dim_H}(A_3) \approx 0.606.
\end{equation*}
Using a ball of functions, $B$, with $\ell_1$ radius $10^{-9}$, centred on a polynomial of degree 120 that we previously proved contains the fixed point, $g$ \cite{Bur20deg3}, after 20 generations of the IFS we obtain rigorous bounds on the Hausdorff dimension of the attractor for the case $d=3$:
\begin{equation}
0.6040883004665372548689708590 < \mathrm{dim_H}(A_3) < 0.6090988814830777819500301340.
\end{equation}

\section{Hausdorff dimension - degree 4 critical point}
\label{sec:d4}

\begin{figure}
\begin{center}
\begin{tabular}{cc}
\includegraphics[width = 5.75cm]{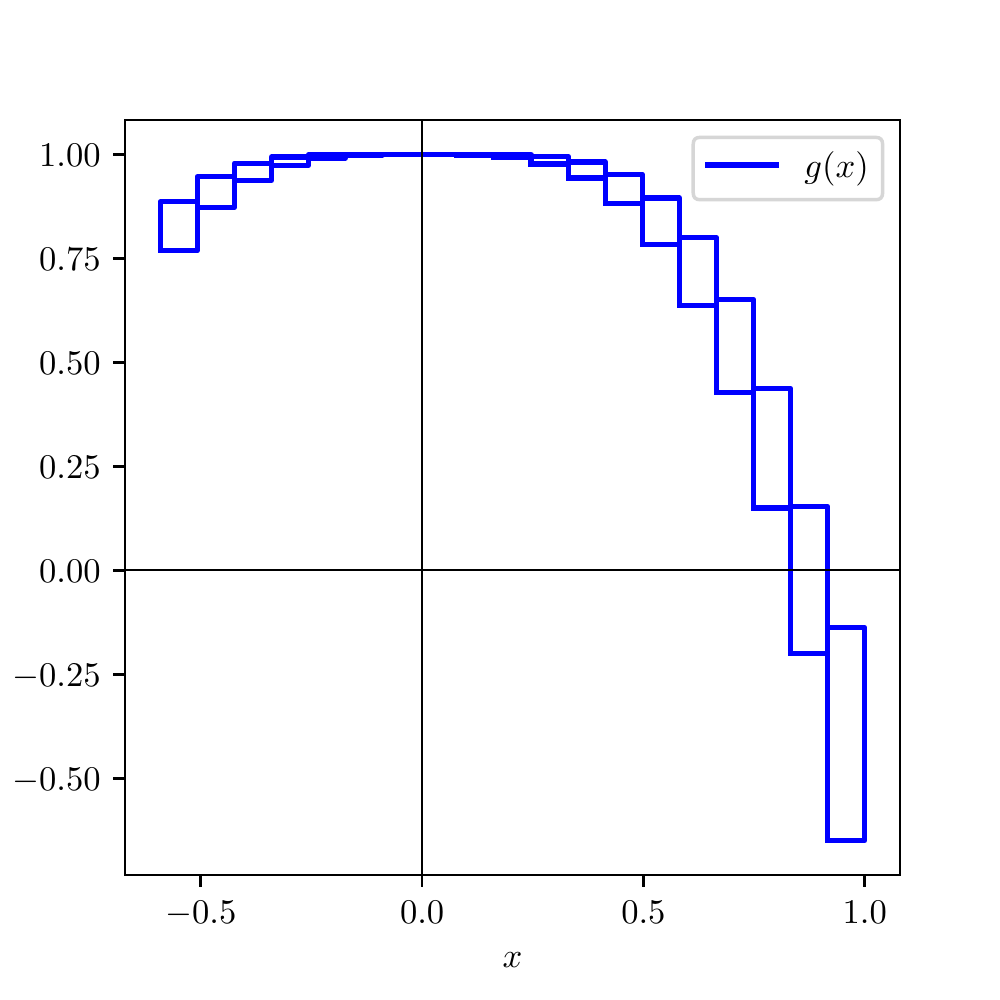} 
& \includegraphics[width = 5.75cm]{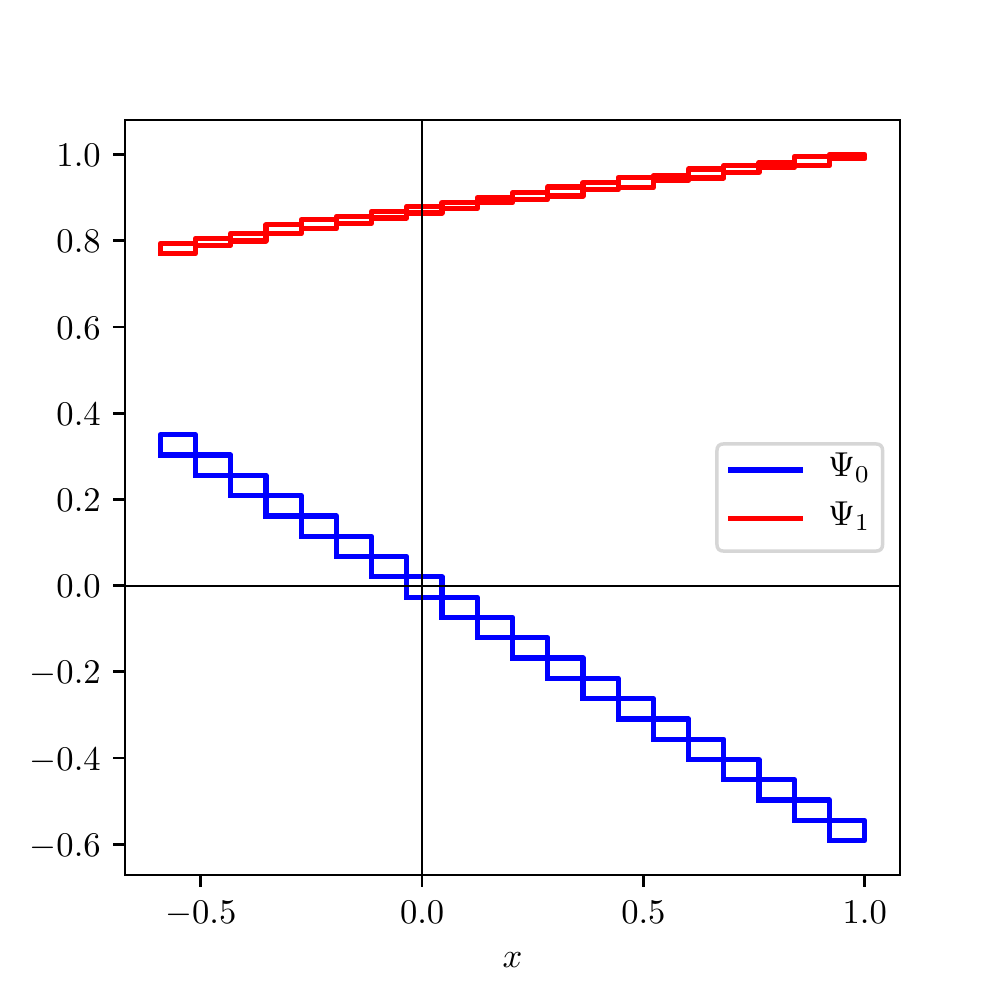}\\
(a) Rigorous bounds on $g$
& (d) Rigorous bounds on $\Psi_0$ and $\Psi_1$\\
    
     \includegraphics[width = 5.75cm]{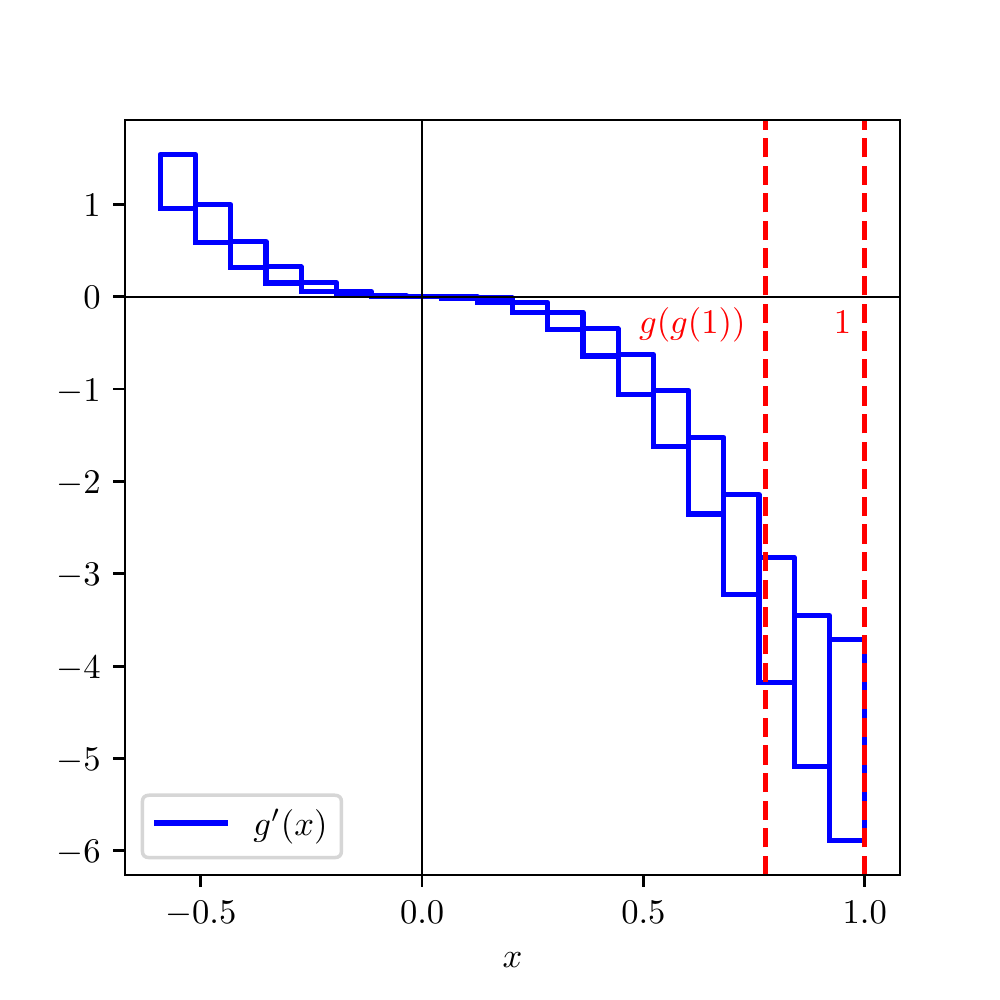} 
     & \includegraphics[width = 5.75cm]{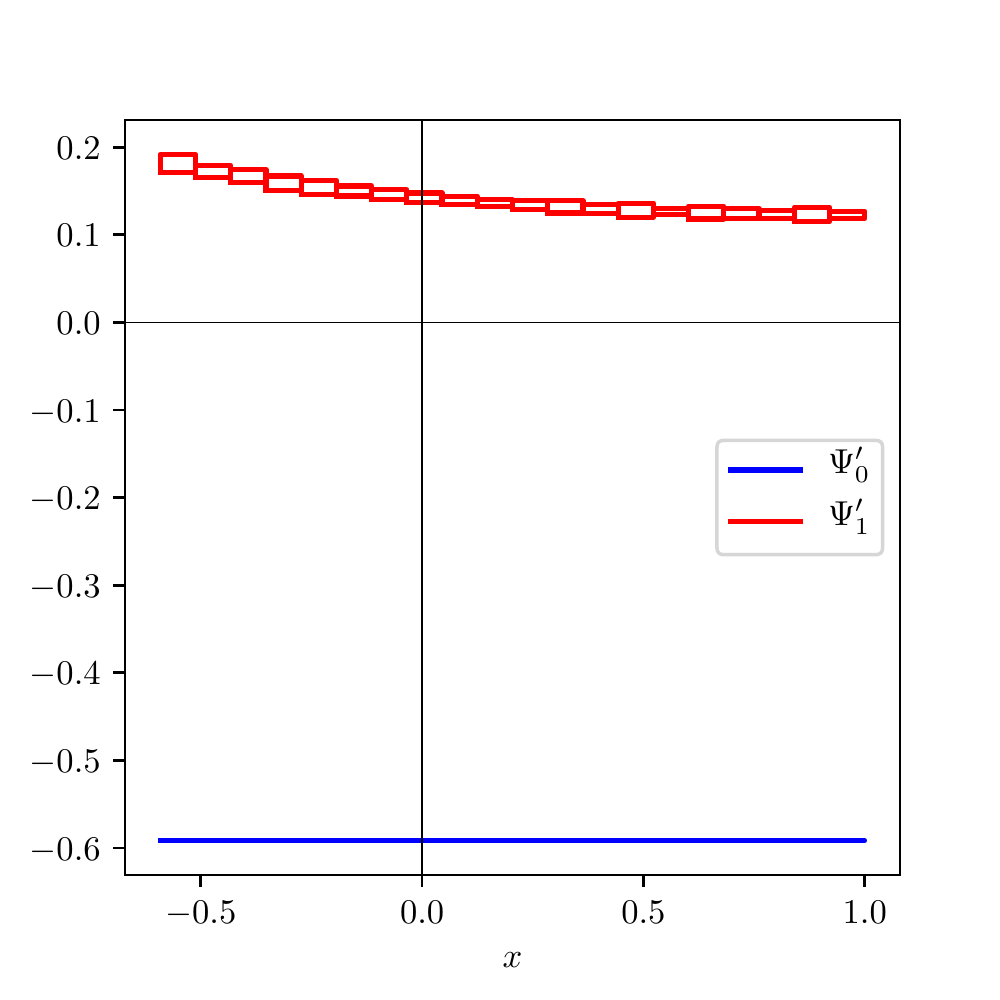}\\
    (b) Rigorous bounds on $g'$
    & (e) Rigorous bounds on $\Psi_0'$ and $\Psi_1'$\\
    
    \includegraphics[width = 5.75cm]{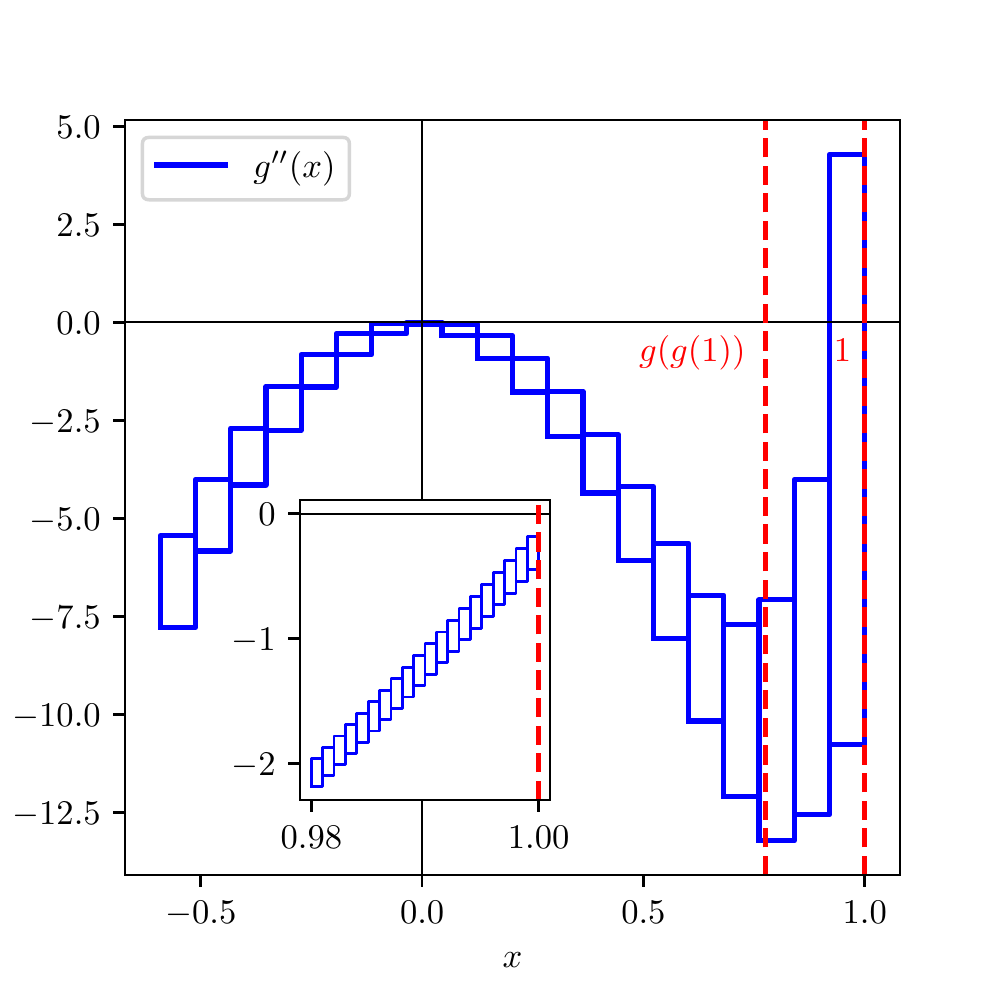} 
    & \includegraphics[width = 5.75cm]{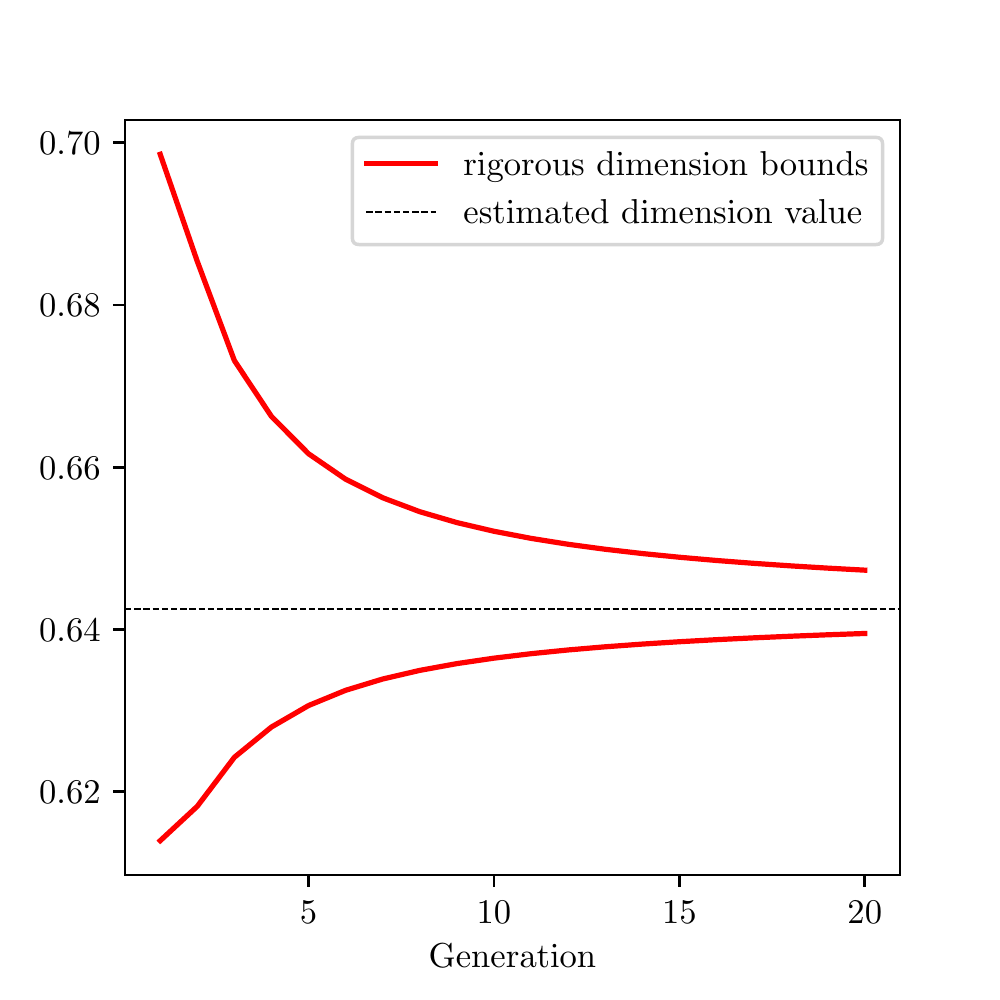}\\
    (c) Rigorous bounds on $g''$
    & (f) Rigorous bounds on dimension
   
\end{tabular}
\caption{Construction of the IFS maps and resulting bounds on dimension for the case $d=4$. The rectangles provide a coarse rigorous covering of the graphs of the corresponding functions: (a) the fixed point function $g$; (b) the first derivative of $g$, with the vertical dashed lines indicating the interval $[g(g(1)),1]$; (c) the second derivative of $g$; from the (coarse) covering shown it is not clear that $g''(x)<0$ on $[g(g(1)),1]$, the inset shows a finer covering of the final segment to establish that $g''(x)<0$; (d) rigorous bounds on  $\Psi_0$ and $\Psi_1$; (e) rigorous bounds on the derivatives of $\Psi_0$ and $\Psi_1$; (f) shows convergence of rigorous upper and lower bounds on dimension with increasing IFS generation.}
\label{Q4Dim}
\end{center}  
\end{figure}

Following the procedure set out for the dimension of attractors of maps with degree 2 critical point, we are able to calculate rigorous bounds on the degree 4 equivalent.

We establish that $g'$ is strictly negative on the interval $[g(g(1)), 1]$ (see figure \ref{Q4Dim}b). However, it is necessary to compute rigorous bounds for $g''$ on much smaller subintervals to establish that $g''$ is strictly negative on the required interval, and therefore that $g'$ is strictly monotonically decreasing (see figure \ref{Q4Dim}c).  The arguments then follow as before to show that $\Psi_1'$ is monotonically decreasing on $[\alpha^{-1},1]$ and we are again able to use the short cut of evaluating derivatives only at the endpoints of the interval.
The bounds produced verify the previously published estimated value for the dimension for maps with a quartic critical point \cite{Kuz02}:
\begin{equation*}
    \mathrm{dim_H}(A_4) \approx 0.642575065.
\end{equation*}
Using a ball of functions, $B$, with $\ell_1$ radius $10^{-9}$, centred on a polynomial of degree 160 that we previously proved contains the fixed point $g$ \cite{Bur20}, after 20 generations of the IFS we obtain rigorous bounds on the Hausdorff dimension of the attractor in the case $d=4$:
\begin{equation}
0.6395131468772885957851829075 < \mathrm{dim_H}(A_4) < 0.6473156929016111051417422795.   
\end{equation}
Figure \ref{Q4Dim}f shows the convergence of the rigorous upper and lower bounds of the Hausdorff dimension.

\section{Conclusion}

Using balls of functions previously proven to contain the renormalisation fixed points for universality classes of maps with degree 2, 3 and 4 critical points, we have calculated rigorous bounds on the Hausdorff dimensions of the corresponding Feigenbaum attractors at the accumulations of period-doubling.

In order to extend this calculation to maps with critical points of degrees higher than 4 a modified version of the method presented above is needed. The short cut described in section \ref{rigFeigDim}, in which we evaluate the derivatives of $\Psi_0$ and $\Psi_1$ only at the endpoints of IFS sub-intervals, cannot be used in the case where the second derivative of $g$ has mixed sign.  Instead, one would either need to partition subintervals into monotonic segments or to bound the suprema and infima on a rigorous covering of the fundamental interval by small sub-intervals. As observed in section \ref{sec:d4}, obtaining tight bounds demands finer subdivisions, increasing the computational cost significantly.

It is also more difficult computationally to bound the renormalisation fixed-point functions themselves as the degree of the critical point is increased: obtaining sufficiently tight bounds requires working with a high truncation degree in a suitable space of analytic functions, while maintaining control over all higher-order terms. The domains of the corresponding power series need to be chosen carefully to ensure that the renormalisation operator is well-defined and differentiable (with compact derivative) on a suitable ball in order that a contraction mapping argument may be used to establish that the ball contains a fixed point.

\end{document}